\documentclass[12pt]{article}%
\usepackage{amsfonts}
\usepackage{sw20bams}
\usepackage{amsmath}
\usepackage{amssymb}
\usepackage{graphicx}%
\providecommand{\U}[1]{\protect\rule{.1in}{.1in}}
\begin{document}

\title{Compositional Square Roots of $\exp(x)$ and $1+x^{2}$}
\author{Steven Finch}
\date{June 7, 2025}
\maketitle

\begin{abstract}
Our work began as an effort to understand calculations by Morris \&\ Szekeres
(1961) and Walker (1991) regarding fractional iteration. \ 

\end{abstract}

\footnotetext{Copyright \copyright \ 2025 by Steven R. Finch. All rights
reserved.}

Let us start by declaring what this paper is not about. \ It is not about
Kneser's 1950 construction \cite{Kn-zeal, BT-zeal} of a real function
$\varphi(x)$, analytic for all $x\in\mathbb{R}$, such that $\varphi
(\varphi(x))=\exp(x)$. \ Among many possible choices of $\varphi(x)$,
Trappmann \&\ Kouznetsov \cite{TK-zeal} found one satisfying a certain
uniqueness criterion and consequently \cite{PC-zeal}%
\[
\varphi(0)=0.49856328794111443467961909249313329400247186492242...,
\]%
\[
\varphi(1)=1.64635423375119458097192403159211451820531164896904...
\]
since $\varphi(1)=\varphi(\varphi(\varphi(0)))=\exp(\varphi(0))$; further
\cite{Cw-zeal, Ps-zeal}%
\[%
\begin{array}
[c]{ccc}%
\varphi(-1)=-0.15588259893582077..., &  & \varphi(-2)=-0.47520567565984653....
\end{array}
\]
A graph of $\varphi(x)$ appears in \cite{TK-zeal}. \ The cost of analyticity
is quite high:\ computations leading to these estimates are elaborate. \ We
turn instead to a different real $\psi(x)$ obeying $\psi(\psi(x))=\exp(x)$
that is less well-known than $\varphi(x)$, due to Szekeres \cite{Sz-zeal}.
\ This choice $\psi(x)$ is infinitely differentiable for all $x\in\mathbb{R}$.
\ Underlying calculations are straightforward and we expand upon the brief
treatment in \cite{MS-zeal}. \ The elegant Mavecha-Laohakosol algorithm
\cite{ML-zeal, F1-zeal}, based on work by de Bruijn \cite{dB-zeal} and
Bencherif \&\ Robin \cite{BR-zeal}, provides greater numerical precision than
was possible in \cite{MS-zeal}. \ Theory due to Walker \cite{W1-zeal}
helpfully clarifies the formulation here. \ We shall also address
half-iterates of $1+x^{2}$, for which a surprisingly simple approach is applicable.

\section{Walker (1991)}

Abel's functional equation%
\[%
\begin{array}
[c]{ccccc}%
g\left(  e^{x}-1\right)  =g(x)+1 &  & \text{becomes} &  & g(y)=g\left(
\ln(1+y)\right)  +1
\end{array}
\]
under the change of variables $y=\exp(x)-1$; hence%
\[
-g\left(  \ln(1+y)\right)  =-g(y)+1.
\]
Consider the recurrence \cite{W1-zeal}%
\[%
\begin{array}
[c]{ccccc}%
y_{n}=\ln\left(  1+y_{n-1}\right)  &  & \text{for }n\geq1, &  & y_{0}>0.
\end{array}
\]
A power-logarithmic series%
\[
y_{n}\sim%
{\displaystyle\sum\limits_{m=0}^{k-1}}
P_{m}\left(  \frac{1}{3}\ln(n)-C\right)  \frac{2}{n^{m+1}}%
\]
is the asymptotic outcome of the Mavecha-Laohakosol algorithm \cite{ML-zeal,
F1-zeal}, valid as $n\rightarrow\infty$, where $P_{m}=P_{m}(X)$ is a
polynomial with rational coefficients:
\[%
\begin{array}
[c]{ccccc}%
P_{0}=1, &  & P_{1}=Y, &  & P_{2}=\dfrac{1}{18}-\dfrac{1}{3}Y+Y^{2},
\end{array}
\]%
\[%
\begin{array}
[c]{ccc}%
P_{3}=-\dfrac{7}{270}+\dfrac{5}{18}Y-\dfrac{5}{6}Y^{2}+Y^{3}, &  &
P_{4}=\dfrac{67}{4860}-\dfrac{53}{270}Y+\dfrac{5}{6}Y^{2}-\dfrac{13}{9}%
Y^{3}+Y^{4},
\end{array}
\]%
\[
P_{5}=-\dfrac{2701}{408240}+\dfrac{653}{4860}Y-\dfrac{83}{108}Y^{2}%
+\dfrac{101}{54}Y^{3}-\dfrac{77}{36}Y^{4}+Y^{5},
\]%
\[
P_{6}=\dfrac{92461}{30618000}-\dfrac{3449}{40824}Y+\dfrac{89}{135}Y^{2}%
-\dfrac{175}{81}Y^{3}+\dfrac{95}{27}Y^{4}-\dfrac{29}{10}Y^{5}+Y^{6}%
\]
and $C=C(y_{0})$ is a constant. \ The parameter $k$ was fixed to be $7$ in
\cite{F1-zeal}; here we fix $k$ to be $20$. \ Our procedure for estimating
$C$, given $y_{0}$, involves computing $y_{N}$ exactly via recursion, for some
suitably large index $N$. We then set the value $y_{N}$ equal to our series
and numerically solve for $C$. \ The assignment $y_{0}\mapsto C(y_{0})$, in
the context of iterations, is written as $y\mapsto-g(y)$ when speaking of
functional equations. \ For example,%
\begin{align*}
-g(1)  &  =2.25696115887251231897468847275855670396572224402525\backslash\\
&  \;\;\;\;\;\;\;71884079911985588153147841049650401141548964254133...
\end{align*}
is accurate to $100$ decimal digits (Walker \cite{W1-zeal} exhibited this
constant to $18$ digits). \ The function $g(y)$ is analytic for all $y>0$
\cite{W2-zeal, W3-zeal}. \ Calculating $g(y)$ presented a \textquotedblleft
major difficulty\textquotedblright\ in \cite{W1-zeal} and we feel very
fortunate to have uncovered the algorithm \cite{ML-zeal} when preparing
\cite{F1-zeal}.

An auxiliary iteration is also defined in \cite{W1-zeal}:%
\[%
\begin{array}
[c]{ccccc}%
h_{n}(x)=\ln\left[  1+h_{n-1}\left(  e^{x}\right)  \right]  &  & \text{for
}n\geq1, &  & h_{0}(x)=x.
\end{array}
\]
The limit%
\[
h(x)=\lim_{n\rightarrow\infty}h_{n}(x)
\]
exists, is infinitely differentiable on $\mathbb{R}$, and satisfies%
\[
h\left(  e^{x}\right)  =e^{h(x)}-1.
\]
Convergence is rapid in this case and \textquotedblleft presents no comparable
difficulty\textquotedblright\ \cite{W1-zeal}; please refer, however, to
Section 7.

\section{Morris \& Szekeres (1961/62)}

From the functional equation \cite{MS-zeal}%
\begin{equation}%
\begin{array}
[c]{ccc}%
A\left(  e^{x}-1\right)  =A(x)+1, &  & A(1)=0
\end{array}
\tag{$\varepsilon$}%
\end{equation}
we deduce that $A(x)=g(x)-g(1)$ and hence%
\begin{align*}
A\left(  h\left(  e^{x}\right)  \right)   &  =g(h\left(  e^{x}\right)
)-g(1)\\
&  =g\left(  e^{h(x)}-1\right)  -g(1)\\
&  =g(h(x))+1-g(1)\\
&  =A(h(x))+1.
\end{align*}
It follows that%
\[
A\left(  h\left(  \exp^{[\sigma+\tau]}(x)\right)  \right)  =A\left(  h\left(
\exp^{[\sigma]}(x)\right)  \right)  +\tau
\]
in general, where $\exp^{[\sigma]}$ denotes the $\sigma^{\text{th}}$ iterate
of $\exp$ and thus%
\[
A\left(  h\left(  \psi(x)\right)  \right)  =A\left(  h\left(  e^{x}\right)
\right)  -\frac{1}{2}%
\]
in particular, where $\psi=\exp^{[1/2]}$ and $\sigma=1$, $\tau=-1/2$. \ For
example, if $x=0$, then%
\begin{align*}
h\left(  e^{0}\right)   &
=1.33030160653615252706823883108246811677165693881095\backslash\\
&  \;\;\;\;\;\ \ 40337638018831101331676581483282140856859692202435...,
\end{align*}%
\begin{align*}
A\left(  h\left(  e^{0}\right)  \right)  -\frac{1}{2}  &
=0.08393456432796234857483098175233421712137805851118\backslash\\
&  \;\;\;\;\;\ \ 97278964441402031609166955346204602065494463863234...,
\end{align*}%
\begin{align*}
h\left(  \psi(0)\right)   &
=1.03762794804982058646342081971330770296371090314246\backslash\\
&  \;\;\;\;\;\ \ 25952307591409488746844256960961601753167791171937...,
\end{align*}%
\begin{align*}
\kappa &  =\psi
(0)=0.49783205633271704965233602445630390782921904427386\backslash\\
&
\;\;\;\;\;\;\;\;\;\;\;\;\;\;\;\ \ \ \ 09031834675751962442930374537608292929637833059997...
\end{align*}
(correcting values $0.0839356$, $1.0376284$, $0.4978330$ in \cite{MS-zeal}).
\ The final result gives immediately%
\begin{align*}
\psi(1)  &  =\psi(\psi(\psi(0)))=\exp\left(  \kappa\right) \\
&  =1.64515080754212070699721442598933345813013919983843\backslash\\
&  \;\;\;\;\;\;\;80498175318207012889797708223177554160246056568719...
\end{align*}
and, because $\lim_{x\rightarrow-\infty}\psi(\psi(x))=0=\psi^{\lbrack
-1]}(\kappa)$,
\begin{align*}
\lim_{x\rightarrow-\infty}\psi(x)  &  =\psi^{\lbrack-2]}(\kappa)=\ln(\kappa)\\
&  =-0.69749249511411060639042389317553029411943210121203\backslash\\
&  \;\;\;\;\;\;\;\;\;\ 80397768357178637152088509859163685260430547444991...
\end{align*}
(correcting values $1.6451523$, $-0.6974906$ in \cite{MS-zeal}). \ Clearly
$\psi(0)\neq\varphi(0)$ and $\psi(1)\neq\varphi(1)$. \ As another example, if
$x=-1$, then%
\begin{align*}
h\left(  e^{-1}\right)   &
=0.97799934154339643974252620984458439993765171689193\backslash\\
&  \;\;\;\;\;\ \ 66327045854679300543898769518123244327448893545092...,
\end{align*}%
\begin{align*}
A\left(  h\left(  e^{-1}\right)  \right)  -\frac{1}{2}  &
=-0.55187537099712606519444010718889813664690793138617\backslash\\
&  \;\;\;\;\;\;\;\;\ \ 32277836533191419967966581021160364680887496829731...,
\end{align*}%
\begin{align*}
h\left(  \psi(-1)\right)   &
=0.80503869740709118874033471320252289938763642733010223\backslash\\
&  \;\;\;\;\;\ \ 43767537970180726819204867556970843716914194961456965...,
\end{align*}%
\begin{align*}
\psi(-1)  &
=-0.15547509757120123423049647987303910430210599991382054\backslash\\
&
\;\;\;\;\;\;\;\;\ \ 83056476810493082392306957428205972651159668227990672....
\end{align*}
As yet another example, if $x=-2$, then%
\begin{align*}
h\left(  e^{-2}\right)   &
=0.88845934112049557219292830881691654769446758907913201\backslash\\
&  \;\;\;\;\;\ \ 90092949534862149927381192238873610051650251323386869...,
\end{align*}%
\begin{align*}
A\left(  h\left(  e^{-2}\right)  \right)  -\frac{1}{2}  &
=-0.78779965475793490231227071926492923554249388942776979\backslash\\
&
\;\;\;\;\;\;\;\;\ \ 42146700660385118863835120598941461877342286411177180...,
\end{align*}%
\begin{align*}
h\left(  \psi(-2)\right)   &
=0.74218536229611431001892860051615014744182602238292869\backslash\\
&  \;\;\;\;\;\ \ 32533345209867321790505797528218757794302473434298452...,
\end{align*}%
\begin{align*}
\psi(-2)  &
=-0.47627409408367607737485701379126035107666779481209780\backslash\\
&
\;\;\;\;\;\;\;\;\ \ 84867647254486020210208099880998544595121523537509950....
\end{align*}
A plot of $\psi(x)$ appears in \cite{MS-zeal}. \ 

\section{A\ Product}

Given $x_{0}=x\in\mathbb{R}$, let $x_{n}=e^{x_{n-1}}$ for all $n\geq1$.
\ Differentiating both sides of%
\[
\exp\left(  h_{n}(x)\right)  =1+h_{n-1}\left(  e^{x}\right)  ,
\]
we obtain%
\[
h_{n}^{\prime}(x)=h_{n-1}^{\prime}\left(  e^{x}\right)  \exp\left[
x-h_{n}(x)\right]
\]
hence%
\[
\frac{h_{n}^{\prime}(x_{0})}{h_{n-1}^{\prime}(x_{1})}=\exp\left[  x_{0}%
-h_{n}(x_{0})\right]  .
\]
More generally,%
\[
\frac{h_{n}^{\prime}(x_{0})}{h_{0}^{\prime}(x_{n})}=\frac{h_{n}^{\prime}%
(x_{0})}{h_{n-1}^{\prime}(x_{1})}\,\frac{h_{n-1}^{\prime}(x_{1})}%
{h_{n-2}^{\prime}(x_{2})}\cdots\frac{h_{2}^{\prime}(x_{n-2})}{h_{1}^{\prime
}(x_{n-1})}\,\frac{h_{1}^{\prime}(x_{n-1})}{h_{0}^{\prime}(x_{n})}=%
{\displaystyle\prod\limits_{j=0}^{n-1}}
\exp\left[  x_{j}-h_{n-j}(x_{j})\right]
\]
and, because $h_{0}^{\prime}(x_{n})=1$,
\[
h^{\prime}(x)=%
{\displaystyle\prod\limits_{j=0}^{\infty}}
\exp\left[  x_{j}-h(x_{j})\right]
\]
in the limit as $n\rightarrow\infty$. \ This infinite product was employed by
Walker \cite{W1-zeal} as a stepping stone to prove that $h(x)$ is infinitely
differentiable. \ We note that%
\[%
\begin{array}
[c]{ccc}%
h(0)=0.8459977040656470559451146..., &  & h^{\prime}%
(0)=0.2893258575528029128912252...
\end{array}
\]
and, as a consequence,%
\[
\left(  g\circ h\right)  ^{\prime}(0)=g^{\prime}\left(  h(0)\right)
h^{\prime}(0)=0.9153668106920727478365342...
\]
which confirms an $8$-digit estimate in \cite{W1-zeal}. \ The product formula
for $h^{\prime}(x)$ will be mentioned again later in Section 6.

\section{Diversions}

In his famous book, Kuczma \cite{Kc-zeal} quoted a result of Baker's
\cite{Ba-zeal} that the equation%
\[
\xi(\xi(x))=\exp(x)-1
\]
possesses no solution which is analytic at $x=0$. \ This theorem is
corroborated by our observation that $A(x)\rightarrow-\infty$ as
$x\rightarrow0^{+}$. Equation ($\varepsilon$) of Section 2 implies
\begin{equation}
A\left(  \xi(x)\right)  =A(x)+\frac{1}{2} \tag{$\delta$}%
\end{equation}
and, since $A(1)=0$,
\[
\xi(1)=A^{[-1]}\left(  \frac{1}{2}\right)  =1.2710274138899515214246518....
\]
(see Table 1 of \cite{MS-zeal}). \ The function $A(x)$ is the unique solution
of ($\varepsilon$) with the property that its derivative $A^{\prime}(x)$ is
\textit{totally monotonic} for every $x>0$ \cite{Sz-zeal}. \ Differentiating
both sides of ($\delta$) yields%
\[
A^{\prime}\left(  \xi(x)\right)  \cdot\xi^{\prime}(x)=A^{\prime}(x)
\]
therefore%
\[
\xi^{\prime}(1)=\frac{A^{\prime}(1)}{A^{\prime}\left(  \xi(1)\right)
}=1.5634166713051621990659414...
\]
employing central difference approximations
\[
A^{\prime}(1)=2.3092256547308712400288373...,
\]%
\[
A^{\prime}\left(  \xi(1)\right)  =1.4770378857499947374615530...
\]
(see Table 2 of \cite{MS-zeal}). \ Also, differentiating both sides of
($\varepsilon$) yields%
\[
A^{\prime}\left(  e^{x}-1\right)  e^{x}=A^{\prime}(x),
\]
a multiplicative analog of ($\varepsilon$). \ For example,%
\[%
\begin{array}
[c]{ccc}%
A\left(  e-1\right)  =A(1)+1=1, &  & A^{\prime}\left(  e-1\right)
=e^{-1}A^{\prime}(1);
\end{array}
\]%
\[%
\begin{array}
[c]{ccc}%
A\left(  \ln(2)\right)  =A(1)-1=-1, &  & A^{\prime}\left(  \ln(2)\right)
=2A^{\prime}(1);
\end{array}
\]%
\[
A\left(  e^{e-1}-1\right)  =A(e-1)+1=2,
\]%
\[
A^{\prime}\left(  e^{e-1}-1\right)  =e^{-(e-1)}A^{\prime}(e-1)=e^{-e}%
A^{\prime}(1);
\]%
\[
A\left(  \ln(1+\ln(2))\right)  =A(\ln(2))-1=-2,
\]%
\[
A^{\prime}\left(  \ln(1+\ln(2))\right)  =(1+\ln(2))A^{\prime}\left(
\ln(2)\right)  =2(1+\ln(2))A^{\prime}\left(  1\right)  .
\]
Finally, the formulas%
\[%
\begin{array}
[c]{ccc}%
\exp\left(  \xi(1)\right)  -1=A^{[-1]}\left(  \dfrac{3}{2}\right)  , &  &
\ln\left(  1+\xi(1)\right)  =A^{[-1]}\left(  -\dfrac{1}{2}\right)
\end{array}
\]
provide a stopping point for us; Paulsen's survey\ \cite{Pa-zeal} will almost
surely give interested readers even greater welcome distraction.

\section{Parable}

The adjective "parabolic" describes both a parable (story containing a
spiritual/ethical lesson) and a parabola (trajectory of a projectile). \ It is
trivially seen that $\left\vert x\right\vert ^{\sqrt{2}}$ is a compositional
square root of $x^{2}$. \ Such a closed-form expression does not exist for the
translate $1+x^{2}$. \ One might expect that an approach involving Abel's
equation is necessary in this scenario. \ In fact, the following simple
iteration suffices \cite{A1-zeal, A2-zeal}:
\[%
\begin{array}
[c]{ccccc}%
f_{n}(x)=\sqrt{f_{n-1}\left(  1+x^{2}\right)  -1} &  & \text{for }n\geq1, &  &
f_{0}(x)=\left\vert x\right\vert ^{\sqrt{2}}.
\end{array}
\]
The limit%
\[
f(x)=\lim_{n\rightarrow\infty}f_{n}(x)
\]
exists and satisfies%
\[
1+f(x)^{2}=f\left(  1+x^{2}\right)  .
\]
This scheme would be doomed if $f_{0}(x)$ failed to effectively approximate
$f(x)$ for large $\left\vert x\right\vert $. \ It reminds us of the recursion
for $h_{n}(x)$ in Section 2. \ One might expect that Walker's \cite{W1-zeal}
theoretical analysis would carry over to here. \ This is again wrong, as will
be seen in Section 6.

For now, our approach will be experimental in flavor. \ The derivative of
$f_{0}(x)$ exists and $f_{0}^{\prime}(0)=0$ in particular; the second
derivative exists apart from the origin:%
\[
\lim\limits_{x\rightarrow0^{-}}f_{0}^{\prime\prime}(x)=+\infty=\lim
\limits_{x\rightarrow0^{+}}f_{0}^{\prime\prime}(x).
\]
The derivative of $f_{1}(x)$ exists apart from the origin:%
\[%
\begin{array}
[c]{ccc}%
\lim\limits_{x\rightarrow0^{-}}f_{1}^{\prime}(x)=-2^{1/4}, &  & \lim
\limits_{x\rightarrow0^{+}}f_{1}^{\prime}(x)=+2^{1/4}.
\end{array}
\]
Beyond this, everything is ideal:\ $f_{n}(x)$ for $n\geq3$ is infinitely
differentiable with vanishing odd derivatives at the origin. \ By uniform
convergence, the same is true for $f(x)$. \ We easily calculate values of $f$
at $0$ and $1$:%
\begin{align*}
f(0)  &  =0.64209450439082838149536305968416471264518420238683\backslash\\
&  \;\;\;\;\;\ \ 22735742257008697493204678006737286460677030322067...,
\end{align*}%
\begin{align*}
f(1)  &  =1.41228535256890352768352266799793732213294564970096\backslash\\
&  \;\;\;\;\;\ \ 98394945437138216556902166907792983626447894164132...;
\end{align*}
its derivative at $1$:%
\begin{align*}
f^{\prime}(1)  &
=1.29305280066488558297825733448798808699131584547046\backslash\\
&  \;\;\;\;\;\ \ 34533926325336590857720303009545158637975557720971...
\end{align*}
and its second/fourth derivatives at $0$:%
\begin{align*}
f^{\prime\prime}(0)  &
=2.01380449734831187813987807106244770668101500994036\backslash\\
&  \;\;\;\;\;\ \ 06389072882388344112810065879668216475636725456000...,
\end{align*}%
\begin{align*}
f^{\prime\prime\prime\prime}(0)  &
=-12.33912516608245045438169032954452268461038223114115\backslash\\
&
\;\;\;\;\;\;\;\;\;\;\;\;90920626627983806171224685669111665774544610174308....
\end{align*}
(the latter via symbolic differentiation of $f_{n}(x)$ for $n>10$). \ 

\section{Two Limits}

The reasoning here is a quadratic variant of Section 3. \ Given $x>0$ WLOG and
$x_{0}=x$, let $x_{n}=1+x_{n-1}^{2}$ for all $n\geq1$. \ Differentiating both
sides of%
\[
f_{n}(x)^{2}=f_{n-1}\left(  1+x^{2}\right)  -1,
\]
we obtain%
\[
2f_{n}(x)f_{n}^{\prime}(x)=f_{n-1}^{\prime}\left(  1+x^{2}\right)  \cdot2x
\]
hence%
\[
\frac{f_{n}^{\prime}(x_{0})}{f_{n-1}^{\prime}(x_{1})}=\frac{x_{0}}{f_{n}%
(x_{0})}.
\]
More generally,%
\[
\frac{f_{n}^{\prime}(x_{0})}{f_{0}^{\prime}(x_{n})}=\frac{f_{n}^{\prime}%
(x_{0})}{f_{n-1}^{\prime}(x_{1})}\,\frac{f_{n-1}^{\prime}(x_{1})}%
{f_{n-2}^{\prime}(x_{2})}\cdots\frac{f_{2}^{\prime}(x_{n-2})}{f_{1}^{\prime
}(x_{n-1})}\,\frac{f_{1}^{\prime}(x_{n-1})}{f_{0}^{\prime}(x_{n})}=%
{\displaystyle\prod\limits_{j=0}^{n-1}}
\frac{x_{j}}{f_{n-j}(x_{j})}%
\]
but, because $f_{0}^{\prime}(x_{n})=\sqrt{2}x_{n}^{\sqrt{2}-1}\neq1$, \
\begin{align*}
f^{\prime}(x)  &  =\lim_{n\rightarrow\infty}\sqrt{2}x_{n}^{\sqrt{2}-1}\cdot%
{\displaystyle\prod\limits_{j=0}^{n-1}}
\dfrac{x_{j}}{f_{n-j}(x_{j})}=\sqrt{2}\lim_{n\rightarrow\infty}\dfrac
{f_{0}(x_{n})}{x_{n}}\cdot%
{\displaystyle\prod\limits_{j=0}^{n-1}}
\dfrac{x_{j}}{f_{n-j}(x_{j})}\\
&  =\sqrt{2}\left(  \lim_{n\rightarrow\infty}\dfrac{1}{x_{n}}\,%
{\displaystyle\prod\limits_{j=0}^{n-1}}
\,x_{j}\right)  \left/  \left(  \lim_{n\rightarrow\infty}\dfrac{1}{f_{0}%
(x_{n})}\,%
{\displaystyle\prod\limits_{j=0}^{n-1}}
\,f_{n-j}(x_{j})\right)  \right.
\end{align*}
cannot be written as a ratio of infinite products (i.e., the representation
for $f^{\prime}(x)$ is less compact than that for $h^{\prime}(x)$). \ Letting%
\[%
\begin{array}
[c]{ccc}%
p(x)=\lim\limits_{n\rightarrow\infty}\dfrac{1}{x_{n}}\,%
{\displaystyle\prod\limits_{j=0}^{n-1}}
\,x_{j}, &  & q(x)=\lim\limits_{n\rightarrow\infty}\dfrac{1}{f_{0}(x_{n})}\,%
{\displaystyle\prod\limits_{j=0}^{n-1}}
\,f_{n-j}(x_{j})
\end{array}
\]
we find
\begin{align*}
p(1)  &  =0.38404642942674944782930679653158458352962978658485\backslash\\
&  \;\;\;\;\;\ \ 27873018833460564662178623304028556617929207228775...,
\end{align*}%
\begin{align*}
q(1)  &  =0.42003208901987412706406902497428715359145947136859\backslash\\
&  \;\;\;\;\;\ \ 80652432802426581656759357241900367167229016820108...
\end{align*}
and verify that $\sqrt{2}p(1)/q(1)$ is indeed $f^{\prime}(1)$ from Section 5
(assessed in a different manner). \ A\ similar verification is possible for
$f^{\prime\prime}(0)$.

\section{Overflow}

In seeking $100$-digit precision for numerical estimates of $h(x)$, we
encountered difficulties that Walker \cite{W1-zeal} could not have
anticipated. \ Return to the exponential setting:\ $x_{n}=e^{x_{n-1}}$ for all
$n\geq1$, where $x_{0}=x$. \ The computer algebra package Mathematica cannot
accept quantities greater than a certain threshold $M$ (called
\textsc{\$MaxNumber}). We observe that $x_{4}<M<x_{5}$ when $x=e^{0}$ or
$x=e^{-1}$, and $x_{5}<M<x_{6}$ when $x=e^{-2}$.\ \ We also see that%
\[%
\begin{array}
[c]{ccc}%
h_{3}\left(  e^{0}\right)  -h_{2}\left(  e^{0}\right)  \approx4.29\times
10^{-9}, &  & h_{4}\left(  e^{0}\right)  -h_{3}\left(  e^{0}\right)
\approx1.84\times10^{-1656529}%
\end{array}
\]
thus $h_{4}\left(  e^{0}\right)  $ offers an excellent approximation of
$h\left(  e^{0}\right)  $, whereas%
\[%
\begin{array}
[c]{ccc}%
h_{3}\left(  e^{-1}\right)  -h_{2}\left(  e^{-1}\right)  \approx
1.03\times10^{-3}, &  & h_{4}\left(  e^{-1}\right)  -h_{3}\left(
e^{-1}\right)  \approx7.08\times10^{-34},
\end{array}
\]%
\[%
\begin{array}
[c]{ccc}%
h_{3}\left(  e^{-2}\right)  -h_{2}\left(  e^{-2}\right)  \approx
4.19\times10^{-3}, &  & h_{4}\left(  e^{-2}\right)  -h_{3}\left(
e^{-2}\right)  \approx3.58\times10^{-13}%
\end{array}
\]
thus $h_{4}\left(  e^{-1}\right)  $, $h_{4}\left(  e^{-2}\right)  $ offer
relatively poor approximations of $h\left(  e^{-1}\right)  $, $h\left(
e^{-2}\right)  $. \ To compute $h_{5}\left(  e^{-1}\right)  $ and
$h_{6}\left(  e^{-2}\right)  $, an asymptotic expansion from \cite{Rg-zeal}:%
\[%
\begin{array}
[c]{ccc}%
\ln(1+x)\sim\ln(x)+\dfrac{1}{x}-\dfrac{1}{2x^{2}}+\dfrac{1}{3x^{3}}, &  &
\text{true as }x\rightarrow\infty,
\end{array}
\]
allows us to safely replace the innermost term $\ln(1+x_{n})$ within
$h_{n}(x)$ by%
\[
\ln(x_{n})+\dfrac{1}{x_{n}}-\dfrac{1}{2x_{n}^{2}}+\dfrac{1}{3x_{n}^{3}%
}=x_{n-1}+\dfrac{1}{x_{n}}-\dfrac{1}{2x_{n}^{2}}+\dfrac{1}{3x_{n}^{3}}%
\]
where $n=5$ or $n=6$. \ All outer terms are left unchanged. \ This replacement
is valid for our purposes since $x_{n}$ is enormous. \ Rather than fatal
`overflow' error flags, Mathematica here gives harmless `underflow' warnings.

\section{Threads}

Return to the logarithmic setting:\ $y_{n}=\ln\left(  1+y_{n-1}\right)  $ for
all $n\geq1$, where $y_{0}=y$. \ One loose thread involves the asymptotic
expansion%
\[
\frac{2}{y_{n}}\sim n-\frac{1}{3}W_{n}+\frac{1}{9n}\left(  W_{n}-\frac{1}%
{2}\right)  +\frac{1}{54n^{2}}\left(  W_{n}^{2}-3W_{n}+\frac{7}{5}\right)
\]
through which $g(y)$ was computed by Walker \cite{W1-zeal, W3-zeal}; the
shorthand $W_{n}$ represents $3g(y)+\ln(n)$. \ Our starting point for proving
this is Mavecha \&\ Laohakosol \cite{ML-zeal}:%
\[
\frac{2}{y_{n}}\sim n+%
{\displaystyle\sum\limits_{m=0}^{k-1}}
T_{m+1}\left(  -\frac{1}{3}\ln(n)+C\right)  \frac{1}{n^{m}}%
\]
where $k=3$ and polynomials \
\[%
\begin{array}
[c]{ccccc}%
T_{1}=Y, &  & T_{2}=-\dfrac{1}{18}-\dfrac{1}{3}Y, &  & T_{3}=\dfrac{7}%
{270}+\dfrac{1}{6}Y+\dfrac{1}{6}Y^{2}%
\end{array}
\]
are actually precursors to $P_{1}$, $P_{2}$, $P_{3}$ mentioned earlier. \ We
obtain%
\begin{align*}
\frac{2}{y_{n}}  &  \sim n+\frac{-\frac{\ln(n)}{3}+C}{1}+\frac{-\frac{1}%
{18}-\frac{1}{3}\left(  -\frac{\ln(n)}{3}+C\right)  }{n}+\frac{\frac{7}%
{270}+\frac{1}{6}\left(  -\frac{\ln(n)}{3}+C\right)  +\frac{1}{6}\left(
-\frac{\ln(n)}{3}+C\right)  ^{2}}{n^{2}}\\
&  =n-\frac{\ln(n)-3C}{3}+\frac{-\frac{1}{2}+\left(  \ln(n)-3C\right)  }%
{9n}+\frac{\frac{7}{5}-3\left(  \ln(n)-3C\right)  +\left(  \ln(n)-3C\right)
^{2}}{54n^{2}}\\
&  =n-\frac{3g(y)+\ln(n)}{3}+\frac{\left(  3g(y)+\ln(n)\right)  -\frac{1}{2}%
}{9n}+\frac{\left(  3g(y)+\ln(n)\right)  ^{2}-3\left(  3g(y)+\ln(n)\right)
+\frac{7}{5}}{54n^{2}}%
\end{align*}
and the result holds. \ Higher accuracy can be found by taking $k=6$ and
utilizing%
\[%
\begin{array}
[c]{ccc}%
T_{4}=-\dfrac{13}{1215}-\dfrac{29}{270}Y-\dfrac{2}{9}Y^{2}-\dfrac{1}{9}%
Y^{3}, &  & T_{5}=\dfrac{305}{81648}+\dfrac{11}{162}Y+\dfrac{127}{540}%
Y^{2}+\dfrac{7}{27}Y^{3}+\dfrac{1}{12}Y^{4},
\end{array}
\]%
\[
T_{6}=-\dfrac{3359}{3402000}-\dfrac{767}{20412}Y-\dfrac{347}{1620}Y^{2}%
-\dfrac{2}{5}Y^{3}-\dfrac{31}{108}Y^{4}-\frac{1}{15}Y^{5}.
\]
We do not know whether $2/y_{n}$ presents advantages over $y_{n}$ when
computing $g(y)$.

The other loose thread surrounds $\psi^{\lbrack-1]}=\ln^{[1/2]}$, which awaits
study. \ Abel's functional equation%
\[%
\begin{array}
[c]{ccccc}%
\tilde{g}\left(  -\ln(1-x)\right)  =\tilde{g}(x)-1 &  & \text{becomes} &  &
\tilde{g}\left(  1-e^{-y}\right)  =\tilde{g}\left(  y\right)  +1
\end{array}
\]
under the change of variables $y=-\ln(1-x)$; we wonder if a parallel line of
reasoning leads to a new $\tilde{\psi}^{[-1]}=\ln^{[1/2]}$, which needn't
necessarily be identical to the old $\psi^{\lbrack-1]}$. \ 

More applications of the Mavecha-Laohakosol algorithm appear in \cite{F2-zeal,
F3-zeal, F4-zeal, F5-zeal}, although its use here in determining compositional
roots seems to be novel.

\section{Addendum}

For the sake of completeness, we include details associated with $x=2$:%
\begin{align*}
h\left(  e^{2}\right)   &
=7.38967388704959756967338001867412264275643003632563\backslash\\
&  \;\;\;\;\;\ \ 77508478967873020720788397090454153674858371882832...,
\end{align*}%
\begin{align*}
A\left(  h\left(  e^{2}\right)  \right)  -\frac{1}{2}  &
=1.78608283926079043089781639609923466374021929085276\backslash\\
&  \;\;\;\;\;\;\;\;84281597367416244327900547868860386215135173504211...,
\end{align*}%
\begin{align*}
h\left(  \psi(2)\right)   &
=3.46490737866885451099147553780807136336994488854327\backslash\\
&  \;\;\;\;\;\ \ 55176242928655544073227449074835977370360415638371...,
\end{align*}%
\begin{align*}
\psi(2)  &  =3.43313194233707999644955791802824671382563886889230\backslash\\
&  \;\;\;\;\;\;\;\;53082295185589703478973854033072734945286302488622....
\end{align*}
as well as tabulated results for $3\leq x\leq10$ and\ $-10\leq x\leq-3$:%
\[%
\begin{tabular}
[c]{|c|c|}\hline
$x$ & $\psi(x)$\\\hline
$3$ & $%
\begin{array}
[c]{c}%
5.99202197295132031000615674624836298613970724943717\backslash\\
\;\;\;\;\;82182963373778442492469192150668270814525104810116...
\end{array}
$\\\hline
$4$ & $%
\begin{array}
[c]{c}%
9.51489757270138263822837109253655450484086775586905\backslash\\
\;\;\;\;\;86894870030500632680086288110521097255427715449911...
\end{array}
$\\\hline
$5$ & $%
\begin{array}
[c]{c}%
14.17243311972469532707908073180616764172011248076733\backslash\\
\;\;\;\;\;\ 21591407609172403468298565131531961558421603058478...
\end{array}
$\\\hline
$6$ & $%
\begin{array}
[c]{c}%
20.13881997410630706401259082693879407757875958428951\backslash\\
\;\;\;\;\;\ 58388999677520912737281561271601677736852854258807...
\end{array}
$\\\hline
$7$ & $%
\begin{array}
[c]{c}%
27.61193319593643272247507619012705816065476880113543\backslash\\
\;\;\;\;\;\ 63435657975625588579603344209001077388104140721944...
\end{array}
$\\\hline
$8$ & $%
\begin{array}
[c]{c}%
36.81357206427293270861578617197111165111529799261040\backslash\\
\;\;\;\;\;\ 12553753149572747145768427516272024247882310757995...
\end{array}
$\\\hline
$9$ & $%
\begin{array}
[c]{c}%
47.98832946439370760636550849427845330058289133011354\backslash\\
\;\;\;\;\;\ 72057598230087486296522641720766204350910277586041...
\end{array}
$\\\hline
$10$ & $%
\begin{array}
[c]{c}%
61.40332280372229980596619288810532944095658957021239\backslash\\
\;\;\;\;\;\ 49914088886531776481803056430249672707570501878582...
\end{array}
$\\\hline
\end{tabular}
\ \ \ \
\]%
\[%
\begin{tabular}
[c]{|c|c|}\hline
$x$ & $\psi(x)$\\\hline
$-3$ & $%
\begin{array}
[c]{c}%
-0.61239998714859786600591623971855133371952340737222\backslash\\
\;\;\;\;\;\ 80141285661438566103305537458505850756855677501818...
\end{array}
$\\\hline
$-4$ & $%
\begin{array}
[c]{c}%
-0.66563384594378810697421323454407784489519860911197\backslash\\
\;\;\;\;\;\ 59202182546222190331385830578021627447074892443190...
\end{array}
$\\\hline
$-5$ & $%
\begin{array}
[c]{c}%
-0.68569425692098230652339835647577505116085769145406\backslash\\
\;\;\;\;\;\ 75314931258837081345515548123168558705279523665760...
\end{array}
$\\\hline
$-6$ & $%
\begin{array}
[c]{c}%
-0.69314143915920659923392521163672797285041339252278\backslash\\
\;\;\;\;\;\ 02599994855473446003115833410264652449740271875808...
\end{array}
$\\\hline
$-7$ & $%
\begin{array}
[c]{c}%
-0.69589037138588905213762541557742644347064918848381\backslash\\
\;\;\;\;\;\ 22097458391710614660301276360157126825760472166030...
\end{array}
$\\\hline
$-8$ & $%
\begin{array}
[c]{c}%
-0.69690290878304236035845043306810790267100968024414\backslash\\
\;\;\;\;\;\ 63126636947978757673809853430773306387327436789226...
\end{array}
$\\\hline
$-9$ & $%
\begin{array}
[c]{c}%
-0.69727557161458021865729407485321647993832157990077\backslash\\
\;\;\;\;\;\ 99613930576876719019917844819952927349323331394021...
\end{array}
$\\\hline
$-10$ & $%
\begin{array}
[c]{c}%
-0.69741268978905283798978761078434988497627641004500\backslash\\
\;\;\;\;\;\ 10846273169142272124441170478664106991576390915960...
\end{array}
$\\\hline
\end{tabular}
\ \ \ \
\]
(correcting the figure $61.3841$ in \cite{MS-zeal} for $x=10$). \ Poetic
(\textquotedblleft special\textquotedblright) values are%
\[%
\begin{tabular}
[c]{|c|c|c|}\hline
$x$ & $\exp^{[1/2]}(x)$ & $\ln^{[1/2]}(x)$\\\hline
$-\infty$ & $\ln(\kappa)$ & .\\\hline
$\ln(\kappa)\approx-0.69749$ & $0$ & $-\infty$\\\hline
$0$ & $\kappa$ & $\ln(\kappa)$\\\hline
$\kappa\approx0.49783$ & $1$ & $0$\\\hline
$1$ & $\exp(\kappa)$ & $\kappa$\\\hline
$\exp(\kappa)\approx1.64515$ & $e$ & $1$\\\hline
$e\approx2.71828$ & $\exp^{[2]}(\kappa)$ & $\exp(\kappa)$\\\hline
\end{tabular}
\
\]
and prosaic values\ $\ln^{[1/2]}(2)$, $\ln^{[1/2]}(3)$ are%
\begin{align*}
\ln(\psi(2))  &
=1.23347294727533348464410093103465477555387622396914\backslash\\
&  \;\;\;\;\;\;\;88852090959674909291949374902295975188119985488451...
\end{align*}%
\begin{align*}
\ln(\psi(3))  &
=1.79042891325614540324620443046860311187798252933100\backslash\\
&  \;\;\;\;\;\;\;15970506143743146039654352404684503005743267773247...
\end{align*}
respectively. \ Clearly $\ln(x)\ll\ln^{[1/2]}(x)\ll x\ll\exp^{[1/2]}(x)\ll
\exp(x)$ as $x\rightarrow\infty$.

\section{Acknowledgements}

I am grateful to Daniel Lichtblau at Wolfram Research for kindly answering my
questions, e.g., about generalizing my original Mathematica code to arbitrary
$k$. \ Michael Rogers \cite{Rg-zeal} solved the problem of evaluating
$h_{5}\left(  e^{0}\right)  $. \ William Paulsen \cite{Ps-zeal, Pa-zeal} and
an\ anonymous contributor \cite{A1-zeal, A2-zeal} assisted me in more ways
than they can imagine. \ The creators of Mathematica earn my gratitude every
day:\ this paper could not have otherwise been written. \ An interactive
computational notebook is available \cite{F6-zeal} which might be helpful to
interested readers.

\end{document}